\documentclass[11pt]{amsart}
\usepackage{amssymb,amsfonts, amsmath, mathrsfs}

\title[Measures of Weak Compactness and Fixed Point Theory]{Measures of Weak Compactness and\\Fixed Point Theory}
\author{Cleon S. Barroso\and Donal O'Regan}
\thanks{The first author is grateful for the financial support by CAPES}
\address{Departamento de Matem\'atica, Universidade Federal do Cear\'a, Campus do Pici, Bl. 914, 60455-760, Fortaleza, CE, Brazil.}
\email{cleonbar@mat.ufc.br}
\address{Department of Mathematics, National University of
Ireland, Galway, Ireland.} \email{donal.oregan@nuigalway.ie}

\keywords{Measure of weak compactness, Fixed point, Sadovskii,
Weakly sequentially continous, Hammerstein itegral equations.}
\date{}

\newtheorem{theorem}{Theorem}[section]

\newtheorem{proposition}{Proposition}[section]

\newtheorem{corollary}{Corollary}[section]
\newtheorem{definition}{Definition}[section]
\numberwithin{equation}{section}

\def\N{\mathbb{N}}
\def\L{\mathscr{L}}

\allowdisplaybreaks

\begin{document}
\subjclass[2000]{Primary 47H10; Secondary 45N05.}

\begin{abstract}
In this paper, we study a class of Banach spaces, called
$\phi$-spaces. In a natural way, we associate a measure of weak
compactness in such spaces and prove an analogue of Sadovskii
fixed point theorem for weakly sequentially continuous maps. A
counter-example is given to justify our requirement. As an
application, we establish an existence result for a Hammerstein
integral equation in a Banach space.
\end{abstract}

\maketitle

\section{Introduction}

Kuratowski measure of noncompactness $\chi (C)$ of a bounded
subest $C$ of a metric space $(X,d)$ is defined to be the infimum
of the set of all $\varepsilon>0$ with the following property:
\begin{center}
C can be covered by finitely many sets, each of whose diameter is
$\leq\varepsilon$.
\end{center}

Let $X$ be a Banach space and let $\mathcal{B}(X)$ denote the
collection of nonempty, bounded subsets of $X$. Measure of
noncompactness appear in various contexts and are very useful
tools in Nonlinear Analysis (see, e.g., \cite{5,6,10,12} and the
references therein). In a general way, under some condition of
condensing, it allows us to obtain fixed points results. The first
result in this direction was Sadovskii's fixed point theorem
\cite{4}. It states that if $M\in\mathcal{B}(X)$ is closed, convex
and $T:M\to M$ is a condensing mapping (that is, $T$ is continuous
and $\chi(T(C))<\chi(C)$ for all $C\in\mathcal{B}(X)$), then $T$
has a fixed point.
\par
Let $\phi\in\mathscr{L}(X)$. The main purpose of this paper is to
investigate what conditions will ensure that the function
$\chi\cdot\phi:\mathcal{B}(X)\to [0,\infty)$ is a measure of weak
compactness ( in the sense that if $C\in\mathcal{B}(X)$ is weakly
closed and $\chi(\phi (C))=0$, then $C$ is weakly compact ) and
then to derive from this fact an analogue of Sadovskii's theorem
for weakly sequentially continuous maps. For the sake of
convenience, $X$ will be called a $\phi$-space if the existence of
such a $\phi$ is verified. A theorem of such type has been
obtained by O'Regan \cite{11}. Nevertheless, the definition of our
measure as well as its use in the study of fixed points for such
maps, seems to be new and original.
\par
This paper is organized as follows: In Section 2, we first
formulate the notion of $\phi$-spaces and give a sufficient
condition for a Banach space to be a $\phi$-space. We will see
that automatically reflexive spaces are $\phi$-spaces. However,
there are nonreflexive $\phi$-spaces (see, Example \ref{ex:2}). In
Section 3, we extend Sadovskii's fixed point theorem to the weak
topology setting. More precisely, we prove that every map weakly
sequentially continuous and $\phi$-condensing (see, Definition
\ref{def:3} ) has a fixed point in nonempty closed, convex and
bounded subsets of a $\phi$-space (see, Theorem \ref{trm:1}).
Also, we will see that $\phi$-condensing is essential for the
conclusion of our main result. Finally, in Section 4, we use our
fixed point theory to establish an existence result for
Hammerstein integral equation in a Banach space.

\section{Preliminaries}
Throughout this paper, $X$ stands for a real Banach space with
norm $\|\cdot\|$. In what follows, we proceed to define the notion
of a $\phi$-space. Let $\mathcal{F}$ be the family of all the
bounded and weakly closed subsets of $X$. If $I$ is the identity
map on $X$, we shall denote by $span\{I\}$ the vector space
generated by $I$.

\begin{definition}\label{def:1} Let $\phi\in\mathscr{L}(X)$ be such that $\phi\not\in span\{I\}$. We will say that $X$ is a $\phi$-space if the following condition is satisfied
\begin{itemize}
\item[$(\mathscr{C})$] if $C\in\mathcal{F}$ and $\phi(C)$ is a
relatively compact, then $C$ is weakly compact.
\end{itemize}
\end{definition}

\paragraph{\it Example 1.}\label{ex:1} Every reflexive space is a $\phi$-space for any $\phi\in\mathscr{L}(X)$ with $\phi\not\in span\{I\}$.
\par
In our next result, we give a sufficient condition for a Banach
space to be a $\phi$-space.

\begin{proposition}\label{prop:1} Let $X$ be a Banach space. If there is a $\phi\in\mathscr{L}(X)$ such that $\phi\not\in span\{I\}$, $Ker(\phi)=\{0\}$ and $\|\phi\sp{-1}y\|\leq c\|y\|$ for all $y\in\phi(X)$, $c>0$, then $X$ is a $\phi$-space.
\end{proposition}

\begin{proof} Take an arbitrary set $C\in\mathcal{F}$ and suppose
that $\phi(C)$ is relatively compact. Let $(C\sb i)$ be a
decreasing sequence of nonempty closed, convex subsets of $C$.
Since each $C\sb i$ is closed and $\|\phi\sp{-1}y\|\leq c\|y\|$ for all $y\in\phi(X)$, we have that each $\phi(C\sb i)$ is closed too.
Thus, $(\phi(C\sb i))$ is a decreasing sequence of nonempty closed
and convex subsets of $\overline{\phi(C)}$. By the finite
intersection property of compact sets we conclude that
$\cap\sb{i=1}\sp{\infty}\phi(C\sb i)\neq\emptyset$. Now since
$\phi\sp{-1}(\phi(C\sb i))=C\sb i$ we get
$\cap\sb{i=1}\sp{\infty}C\sb i\neq\emptyset$. Applying now the
$\check{S}mulian$ nested interval principle as in \cite[pp.
433]{2}, one has that $C$ is weakly compact.
\end{proof}

\paragraph{\it Remark.} From the open mapping theorem we know that if $\phi\in\mathscr{L}(X)$ is one-to-one and $\phi(X)$ is a Banach space then there exists $\,c>0\,$ with $\,\|\phi\sp{-1}y\|\leq c\|y\|\,$ for $\,y\in X$.

\paragraph{\it Example 2.}\label{ex:2} This is an example of a $\phi$-space that is not reflexive. Let $C[0,1]$ be the space of the real-valued continuous functions and fix $m\in C[0,1]$ a increasing continuous function such that $m(t)\geq a$, for all $t\in [0,1]$ with $0<a<1$. Let us define $\phi\in\mathscr{L}(C[0,1])$ by $\phi(x)(t)=m(t)\cdot x(t)$. Then $\phi\not\in span\{I\}$ and $Ker(\phi)={0}$. Moreover, if $(Sy)(t)=m(t)\sp{-1}\cdot y(t)$ for all $y\in C[0,1]$, then $\|Sy\|\sb{\infty}\leq\|1/m\|\sb{\infty}\cdot \|y\|\sb{\infty}$ for every $y\in C[0,1]$ and $S\sb{|\phi(C[0,1])}\equiv \phi\sp{-1}$. By Proposition \ref{prop:1}, $C[0,1]$ is a $\phi$-space.
\par
We are now going to define a natural measure of weak  compactness
in a $\phi$-space. For convenience, we first summarize the main
properties of $\chi$ which will be used here.

The function $\chi :\mathcal{B}(X)\to [0,1]$ enjoys the following
properties:
\begin{itemize}
\item[(1)] For any $C\in\mathcal{B}(X)$, $0\leq \chi (C)\leq\delta
(C)=$diameter of $C$, \item[(2)] $\chi (C)=0$ iff $C$ is
relatively compact, \item[(3)] if $C\subset D$, then
$\chi(C)\leq\chi(D)$ for all $C,D\in\mathcal{B}(X)$,
\item[(4)] $\chi(C\cup D)=\max\{\chi(C),\chi(D)\}$, $\forall
C,D\in\mathcal{B}(X)$, \item[(5)] $\chi(\overline{C})=\chi(C)$,
\item[(6)] $\chi(\overline{co}(D))=\chi(co(C))=\chi(C)$.
\end{itemize}
The details of $\chi$ and their properties may be found in
\cite{1}.
\par

\begin{definition}\label{def:2} Let $X$ be a $\phi$-space. Define the measure of weak compactness on $X$ as being the function $\chi\sb\phi:\mathcal{B}(X)\to [0,\infty)$ given by
\begin{center}
$\chi\sb\phi(C)=\chi(\phi(C))$.
\end{center}
\end{definition}

\paragraph{\it Remark.} Clearly, $\chi\sb\phi$ is well defined. Also, it is
follows from the definition of $\phi$-spaces and by $(2)$ that a
set $C\subset\mathcal{F}$ is weakly compact if
$\chi\sb\phi(C)=0$. In addition, we have
$\chi\sb\phi(C)\leq\|\phi\|\chi(C)$ for all $C\in\mathcal{B}(X)$.
Thus if $C$ is compact then $\chi\sb\phi(C)=0$. Hence, in
$\phi$-spaces every compact and weakly closed set is weakly
compact.

\section{Fixed Point Theory}
Before going to the main result of this paper, we give a useful
definition. In the proof of Sadovskii theorem,
$\chi(T(C))<\chi(C)$ plays a crucial role.
\par
This property suggests the following definition.
\begin{definition}\label{def:3} Let $M$ be a bounded subsest of a $\phi$-space $X$. We will say that $T:M\to M$ is ws $\phi$-condensing (resp. w $\phi$-condensing ) if $T$ is weakly sequentially continuous (resp. weakly continuous ) and $\phi$-condensing, i.e,
\begin{center}
$\chi\sb\phi(T(C))<\chi\sb\phi(C)$\quad for  all\quad
$C\in\mathcal{B}(M)$,
\end{center}
where $\mathcal{B}(M)$ denote the collection of all nonempty and
bounded subsets of $M$.
\end{definition}

Our main result is as follows.
\begin{theorem}\label{trm:1} Let $X$ be a $\phi$-space and $M$ a closed, convex and bounded subset of $X$. If $T:M\to M$ is a ws $\phi$-condensing mapping, then there is a $y\in M$ such that $T(y)=y$.
\end{theorem}

Since every weakly continuous mapping is weakly sequentially
continuous we immediately obtain the following corollary.
\begin{corollary}\label{cor:1} Assume $M$ is as in Theorem \ref{trm:1}. If $T$ is w $\phi$-condensing from $M$ into itself, then $T$ has a fixed point.
\end{corollary}

The proof of Theorem \ref{trm:1} is based in the following
generalization of the Schauder-Tychonov fixed point principle
which was obtained by Arino, Gautier and Penot \cite{7}.

\begin{theorem}\label{trm:2} Let $M$ be a nonempty, convex and weakly compact subset of a Banach space $X$ and let $T:M\to M$ be a weakly sequentially continuous operator. Then $T$ has at least one fixed point in the set $M$.
\end{theorem}

\paragraph{\it Proof of Theorem \ref{trm:1}.} We proceed as in the proof of the
Sadovskii's Theorem. Fix a point $x\in M$ and let $\Sigma$ denote
the system of all closed, convex subsets of $M$ for which $x\in K$
and $T(K)\subset K$. Now set
\begin{center}
$B=\cap\sb{K\in\Sigma}K$, \qquad
$C=\overline{co}\{T(B)\cup\{x\}\}$.
\end{center}
Then $B=C$, $T(C)\subset C$ and
$\phi(C)\subset\overline{co}\{\phi(T(B))\cup\{\phi(x)\}\}$. Using
$(3)$,$(4)$ and $(6)$, one has
\begin{center}
$\chi(\phi(C))\leq\chi(\phi(T(B)))=\chi(\phi(T(C))).$
\end{center}
Since $T$ is $\phi$-condensing, it follows that
$\chi\sb\phi(C)=0$. Then in view of the definition of
$\phi$-spaces we conclude that $C$ is weakly compact since
$C\in\mathcal{F}$. Now Theorem \ref{trm:2} gives a fixed point
$y\in C$ for $T$. This ends the proof.\hspace{9cm} $\square$

Evidently, if $X$ is a reflexive space the assumption
$\chi\sb\phi(T(C))<\chi\sb\phi(C)$ in the Theorem \ref{trm:1} is
unnecessary. But, in a nonreflexive $\phi$-space this requirement
may be essential as the next example shows.

\paragraph{\it Example 4.} Consider the Banach space $X=l\sp 1$ endowed with its usual norm. Let us define $\phi: l\sp 1\to l\sp 1$ by 
$$
\phi\big( \sum\sb{n\in\N}\alpha\sb n e\sb n\big )=\sum\sb{n\in\N}\alpha\sb n e\sb{n+1}.
$$
Then $\phi\in\L(l\sp 1)$, $Ker(\phi)=\{0\}$ and $\phi\not\in span\{I\}$. Now, define $S:l\sp 1\to l\sp 1$ by
$$
S\big(\sum\sb{n\in\mathbb{N}}\alpha\sb n e\sb n\big)=\sum\sb{n\geq 2}\alpha\sb n e\sb{n-1}.
$$ Then, $S\sb{|\phi(l\sp 1)}\equiv \phi\sp{-1}$ and $\|S(\eta)\|\leq \|\eta\|$ for all $\eta\in l\sp 1$. Consequently, $\|\phi\sp{-1}y\|\leq \|y\|$ for all $y\in \phi(l\sp 1)$. By Propostion \ref{prop:1}, $l\sp 1$ is a $\phi$-space. Let now $M$ be the set
$$
M=\big\{\sum\sb{n\in\N}\alpha\sb n e\sb n : \quad \alpha\sb n\geq 0, \quad \sum\sb{n\in\N} \alpha\sb n=1\big\}.
$$
Define the map $T:M\to M$ as $T=\phi\sb{| M}$. Then, $T$ is fixed point free on $M$. Moreover, one can check that $T$ is an isometry on $M$. Also it is easy to see that $T$ is weakly sequentially continous, and 
$$
\|\phi T(\eta)-\phi T(\mu)\|=\|\phi(\eta)- \phi(\mu)\|,
$$
for all $\mu, \eta \in M$, which implies that $\chi\sb\phi(TC)\geq \chi\sb\phi(C)$ for all bounded set $C$ in $M$.

\paragraph{\it Remark.} The counterexample above was inspired by the paper \cite[pp. 444]{3}.

\section{Hammerstein integral equations in Banach spaces}
Let $E$ be a Banach space and consider the Hammerstein integral
equation
\begin{equation}\label{eqn:1}
y(t)=h(t)+\int\sb 0\sp 1 k(t,s)f(s,y(s))ds,\quad t\in [0,1],
\end{equation}
where
\begin{equation}\label{eqn:2}
\left\{
\begin{split}
&f:[0,1]\times E\to E \mbox{ is such that, } f\sb t=f(t,\cdot):E\to E\mbox{ is} \\
&\mbox{weakly sequentially continuous, for each } t\in  [0,1],
\end{split}\right.
\end{equation}

\begin{equation}\label{eqn:3}
h\in C([0,1],E)\mbox{ is arbitrary },
\end{equation}
and
\begin{equation}\label{eqn:4}
\left\{
\begin{split}
& k(t,s)\in L\sp 1([0,1],\mathbb{R}) \mbox{ for each } t\in [0,1] \mbox{ and}\\
& \mbox{the map }t\mapsto k(t,s) \mbox{ is continuous from }[0,1]
\mbox{ to } L\sp 1([0,1],\mathbb{R})
\end{split}\right.
\end{equation}
hold.
\par
When $E$ is a reflexive space and $f$ is weakly-weakly continuous,
O'Regan \cite{9} proved that $(\ref{eqn:1})$ has a solution if the
additional conditions
\begin{equation}\label{eqn:5}
\left\{
\begin{split}
&\mbox{there exists a nondecreasing continuous function}\\
&\Omega:[0,\infty)\to (0,\infty)\mbox{ with }
\|f(s,x)\|\leq\Omega(\|x\|)\mbox{ for } t\in[0,1]
\end{split}\right.
\end{equation}
and
\begin{equation}\label{eqn:6}
\big( K=\sup\sb{t\in[0,1]}\int\sb 0\sp 1 |k(t,s)|ds\big )\cdot
\limsup\sb{t\to\infty}\frac{\Omega(t)}{t}\leq 1,
\end{equation}
are satisfied. The proof in \cite{9} depends strongly on the
reflexivity of $E$ since a weak version of the Arzela-Ascoli
theorem is used to guarantee that the right-side in (\ref{eqn:1})
defines a weakly compact operator on $C([0,1],E)$. However, in an
arbitrary Banach space $E$ this technique may fails due to lack of
weak compactness in balls. Thus it is natural to ask what
conditions are needed if the word reflexive is removed. Here as an
application of our fixed point theory, we offer the following
result as a partial answer to this question.

\begin{theorem}\label{trm:3} Let $E$ be a Banach space and assume $(\ref{eqn:2})$-$(\ref{eqn:6})$. In addition suppose that for some $\phi\in\mathscr{L}(E)$ with $\phi\not\in span\{I\}$, $Ker(\phi)=\{0\}$ and $\|\phi\sp{-1}(y)\|\leq c\|y\|$ for all $y\in\phi(E)$ with $c>0$ we have
\begin{equation}\label{eqn:7}
\|\phi(\int\sb 0\sp 1 k(t,s)(f(s,x)-f(s,y))ds)\|\sb
E\leq\alpha\|\phi(x(t)-y(t))\|\sb E,
\end{equation}
for all $t\in[0,1]$, $x,y\in C([0,1],E)$ and some $0<\alpha<1$.
Then, $(\ref{eqn:1})$ has a solution $y\in C([0,1],E)$.
\end{theorem}

\begin{proof} Let $X=C([0,1],E)$ with its usual norm $\|\cdot\|$.
Let us define $\phi\in\mathscr{L}(X)$ by $(\phi x)(t)=\phi(x(t))$.
Then, $\phi\not\in span\{I\}$, $Ker(\phi)=\{0\}$ and $\|\phi\sp{-1}y\|\leq c\|y\|$ for all $y\in\phi(C(I,E))$. Proposition
\ref{prop:1} now implies that $X$ is a $\phi$-space. Now
$(\ref{eqn:5})$,$(\ref{eqn:6})$ guarantees that there exists $R>0$
such that
\begin{equation}\label{eqn:8}
\|h\|+K\Omega(R)\leq R,
\end{equation}
(cf. \cite{9} ). Let $M$ be the set
\begin{center}
$M=\{x\in X:\quad \|x\|\leq R\}.$
\end{center}
Let us define the map $T:M\to M$ by
$$
(Tx)(t)=h(t)+\int\sb 0\sp 1 k(t,s)f(s,x(s))ds,\quad t\in [0,1].
$$

We claim now that $T$ is weakly sequentially continuous. Indeed,
let $(x\sb n)\subseteq M$ such that $x\sb n\rightharpoonup x$ in
$X$. Then, $x\sb n(s)\rightharpoonup x(s)$ in $E$ for each
$s\in[0,1]$. By (\ref{eqn:2}), one has that $f(s,x\sb
n(s))\rightharpoonup f(s,x(s))$ in $E$ for each $s\in[0,1]$, as
$n\to\infty$. Combining now the Hahn-Banach theorem with the
Dominated convergence theorem we conclude that $(Tx\sb n)(t)\to
(Tx)(t)$ in $E$ for each $t\in[0,1]$, as $n\to\infty$. Now, from
(\ref{eqn:4}) it is easy to check that $\{ (Tx\sb
n):~n\in\mathbb{N}\}$ is an equicontinuous subset of $X$. Thus, by
Ascoli-Arzela theorem, we may conclude that
$Tx\sb{n\sb j}\to Tx$ in $X$, for some subsequence $(x\sb{n\sb j})$ of $(x\sb n)$. Hence $Tx\sb n\rightharpoonup Tx$ in
$X$, as $n\to\infty$. If not there exists a subsequence $(Tx\sb{n\sb k})$ of $(Tx\sb n)$ and a weak neighborhood $V\sp w(Tx)$ of $Tx$ such that $Tx\sb{n\sb k}\not\in V\sp w(Tx)$, for all $k\in\mathbb{N}$. Now arguing as before we find a subsequence $Tx\sb{n\sb{k\sb l}}\rightharpoonup Tx$, which is a contradiction. This proves the claim. Further, in view of
(\ref{eqn:7}) it follows that $T$ is $\phi$-condensing. Now, in
order to apply Theorem \ref{trm:1}, it remains to show that
$T(M)\subset M$. But (\ref{eqn:8}) implies that $\|Tx\|\leq
\|h\|+K\Omega(R)\leq R$ whenever $x\in M$. Also, from
(\ref{eqn:3}) and (\ref{eqn:4}) one easily verifies that $Tx$ is
norm continuous if $x\in M$. This completes the proof of the
Theorem.
\end{proof}

\paragraph{\it Remark.} Notice (\ref{eqn:6}) could be replaced by condition:
there exists $R>0$ with $\|h\|+K\Omega(R)\leq R$.


\begin{thebibliography}{10}
\bibitem{1}
E. Zeidler,
\newblock{\em Nonlinear Functional Analysis and its Applications: Fixed point Theorems,}
\newblock{Vol. I, Spinger Verlag, New York, Berlin, Heidelberg, Tokyo}, (1986). MR {\bf 87f:} 47083

\bibitem{2}
Dunford, N. and Schwartz, J.
\newblock{\em Linear Operators,}
\newblock{Vol. I, Interscience, New York,} (1958).

\bibitem{3}
P.N. Dowling and C.J. Lennard,
\newblock{\em Every noreflexive subspace of $L\sp 1[0,1]$ fails the fixed point property,}
\newblock{Proc. Amer. Math. Soc.,}  {\bf 125} (1997), 443--446. MR {\bf 97d:} 46034

\bibitem{4}
B. N. Sadovskii, 
\newblock{\em A fixed point principle,}
\newblock{Functional Analysis and its Applications,} {\bf 1} (2), (1967), 151--153. MR {\bf 35} \#2184

\bibitem{5}
J. Bana\'s,
\newblock{\em Applications of measure of weak noncompactness and some classes of operators in the theory of functional equations in the Lebesgue space, Proceedings of the Second World Congress of Nonlinear Analysts,}
\newblock{Part 6 (Athens, 1996), Nonlinear Anal.,} {\bf 30} (1997), 3283--3293. MR {\bf 99b:} 47091

\bibitem{6}
J. Bana\'s and A. Martin\'on,
\newblock{\em Measures of weak noncompactness in Banach sequence spaces,}
\newblock{Portugal. Math.,} {\bf 52} (1995), 131--138. MR {\bf 96i:} 46010

\bibitem{7}
O. Arino, S. Gautier, and J.P. Penot, 
\newblock{\em A fixed point theorem for sequentially continuous mappings with applications to ordinary differential equations,}
\newblock{Funkc. Ekvac.,} {\bf 27} (1984), 273--279. MR {\bf 86g:} 47069

\bibitem{8}
Reed, M. and Simon, B.,
\newblock{\em Methods of Modern Mathematical Physics,}
\newblock{Academic, New York,} Vol.I (1971). MR {\bf 58} \#12429a

\bibitem{9}
D. O'Regan, 
\newblock{\em Integral equations in reflexive Banach spaces and weak topologies,}
\newblock{Proc. Amer. Math. Soc.,} {\bf 124} (1996), 607--614. MR {\bf 96d:} 45003

\bibitem{10}
D. O'Regan, 
\newblock{\em Fixed-Point Theory for the Sum of Two Operators,}
\newblock{Appl. Math. Lett.,} {\bf 9} no.1 (1996), 1--8. MR {\bf 97e:} 47101

\bibitem{11}
D. O'Regan, 
\newblock{\em Fixed point theory for weakly sequentially continuous mappings,}
\newblock{Math. Comput. Modelling,} {\bf 27} (5), (1998), 1--14. MR {\bf 99c:} 47092

\bibitem{12}
J.M. Ayerbe Toledano, T. Domm\'\i nguez Benavides and G. L\'opez
Acedo, 
\newblock{\em Measures of noncompactness in metric fixed point theory,}
\newblock{Birkh\"auser Verlag, Basel Boston Berlin,} 1997. MR {\bf 99e:} 47070
\end{thebibliography}
\end{document}